\documentstyle[12pt,amssymb,graphicx]{article}

\begin{document}
\renewcommand{\thefootnote}{}
\date{}

\def\thebibliography#1{\noindent{\normalsize\bf References}
 \list{{\bf
 \arabic{enumi}}.}{\settowidth\labelwidth{[#1]}\leftmargin\labelwidth
 \advance\leftmargin\labelsep
 \usecounter{enumi}}
 \def\newblock{\hskip .11em plus .33em minus .07em}
 \sloppy\clubpenalty4000\widowpenalty4000
 \sfcode`\.=1000\relax}

\baselineskip=18pt
\title{
{
\normalsize
{$C_k$-MOVES ON SPATIAL THETA-CURVES AND VASSILIEV INVARIANTS}}}

\author{{\normalsize AKIRA YASUHARA}\\
{\small Department of Mathematics, Tokyo Gakugei University}\\[-1mm]
{\small Nukuikita 4-1-1, Koganei, Tokyo 184-8501, Japan}\\
{\small {\em Current address}, October 1, 1999 to September 30, 2001:}\\[-1mm]
{\small Department of Mathematics, The George Washington University}\\[-1mm]
{\small Washington, DC 20052, USA}\\[-1mm]
{\small e-mail: yasuhara@u-gakugei.ac.jp}}

\maketitle

\vspace*{-5mm}  
\baselineskip=15pt
{\small 
\begin{quote}
\begin{center}A{\sc bstract}\end{center}
\hspace*{1em} The $C_k$-equivalence is an equivalence relation 
generated by $C_k$-moves defined by Habiro. 
Habiro showed that the set of $C_k$-equivalence classes 
of the knots forms an abelian group under the connected sum and 
it can be classified by the additive Vassiliev invariant of order 
$\leq k-1$. We see that the set of $C_k$-equivalence classes 
of the spatial $\theta$-curves forms a group under the vertex connected 
sum and that if the group is abelian, then it can be classified by the 
additive Vassiliev invariant of order $\leq k-1$. 
However the group is not necessarily abelian. In fact, we show that 
it is nonabelian for $k\geq 12$. 
As an easy consequence, we have the set of $C_k$-equivalence classes
of $m$-string links, which forms a group under the composition, 
is nonabelian for $k\geq 12$ and $m\geq 2$. 
\end{quote}}


\renewcommand{\thefootnote}{*}

\noindent
{\bf 1. $C_k$-moves and Vassiliev invariants of spatial 
$\theta$-curves}

\medskip
A {\it tangle} $T$ is a disjoint union of properly embedded 
arcs in the unit $3$-ball $B^{3}$. 
A {\it local move} is a pair of 
tangles $(T_{1},T_{2})$ with $\partial T_{1}=\partial T_{2}$ 
such that for each component $t$ of $T_1$ there exists 
a component $u$ of $T_2$ with $\partial t=\partial u$. 
Two local moves $(T_{1},T_{2})$ and $(U_{1},U_{2})$
are {\it equivalent} 
if there is an orientation preserving 
self-homeomorphism $\psi :B^{3}\rightarrow B^{3}$ such that $\psi (T_{i})$ 
and $U_{i}$ are ambient isotopic in $B^3$ relative to $\partial
B^{3}$ for 
$i=1,2$. Here $\psi (T_{i})$ 
and $U_{i}$ are {\it ambient isotopic in $B^3$ relative to $\partial
B^{3}$} if $\psi (T_{i})$ is deformed to $U_{i}$ by an ambient isotopy of
$B^3$ that is
pointwisely fixed on $\partial B^3$.

Let $(T_{1},T_{2})$ be a 
local move, $t_{1}$ a component of $T_1$ and $t_2$ a component of $T_2$ 
with $\partial t_{1}=\partial t_{2}$. Let $N_1$ and $N_2$ be regular 
neighbourhoods of $t_1$ and $t_2$ in $(B^3-T_1)\cup t_1$ and $(B^3-T_2)\cup
t_2$ respectively such
that $N_{1}\cap
\partial  B^{3}=N_{2}\cap \partial B^{3}$. Let $\alpha$ be a disjoint union
of properly 
embedded arcs in $B^{2}\times [0,1]$ as illustrated in Fig. 1.1.
Let $\psi_{i}:B^{2}\times [0,1]\rightarrow N_{i}$ be a homeomorphism 
with $\psi_{i}(B^{2}\times \{ 0,1\} )=N_{i}\cap \partial B^{3}$ for $i=1,2$. 
Suppose that $\psi_{1}(\partial \alpha )=\psi_{2}(\partial \alpha )$ and 
$\psi_{1}(\alpha )$ and $\psi_{2}(\alpha )$ are ambient isotopic in $B^{3}$ 
relative to $\partial B^3$. Then we say that a local move 
$((T_{1}-t_{1})\cup \psi_{1}(\alpha ), (T_{2}-t_{2})\cup \psi_{2}(\alpha ))$ 
is a {\it double} of $(T_{1},T_{2})$ 
with respect to the components $t_1$ and $t_2$. 

\begin{center} 
\includegraphics[trim=0mm 0mm 0mm 0mm, width=.25\linewidth]
{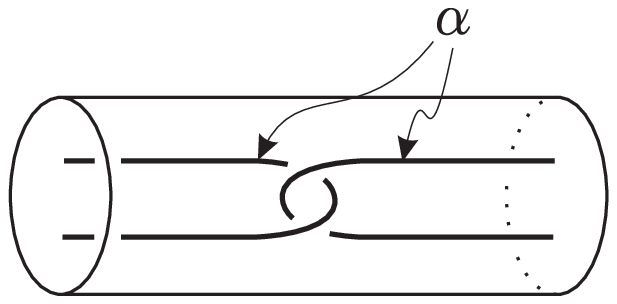}

Fig. 1.1
\end{center}

A {\it $C_{1}$-move} is a local move $(T_1,T_2)$ as illustrated in Fig. 1.2. 
A double of a $C_{k}$-move is called a {\it $C_{k+1}$-move}. 
Note that, for each natural number $k$, 
there are only finitely many $C_{k}$-moves up to 
equivalence. 
It is easy to see that if $(T_1,T_2)$ is a $C_n$-move, then 
$(T_2,T_1)$ is equivalent 
to a $C_{n}$-move (but possibly not equivalent to itself).
The definition of $C_k$-move follows that in \cite{Habiro2}, 
and is defferent from the one in \cite{Habiro1}. 
However by an easy induction on $k$ it is shown that these 
two definitions are essentially same. In \cite{Habiro1}, 
a $C_k$-move is called a {\it simple} $C_k$-move, and 
a $C_k$-move means a {\it parallel} of a $C_k$-move. 
The definition of parallel of a local move appears in Section 3. 

\begin{center} 
\includegraphics[trim=0mm 0mm 0mm 0mm, width=.35\linewidth]
{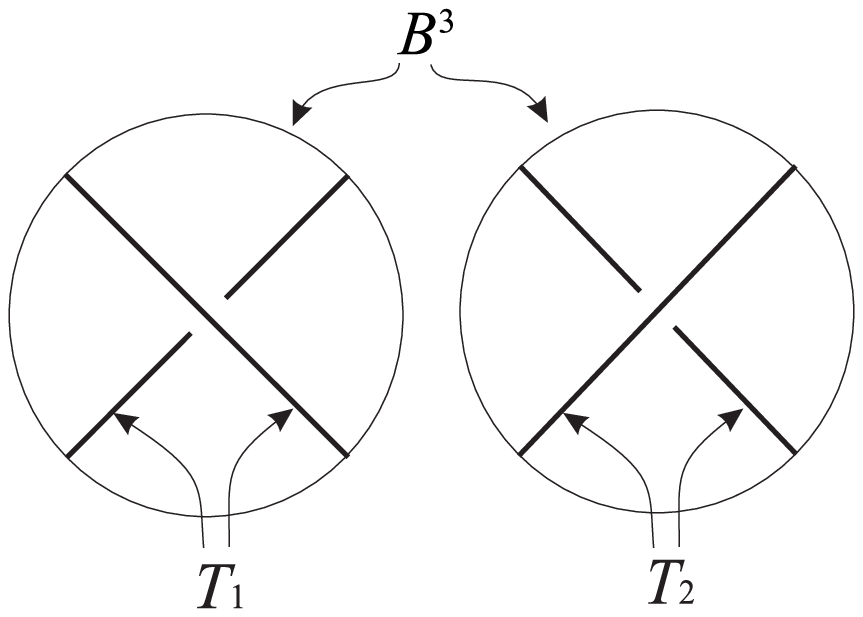}

Fig. 1.2
\end{center}

Let $G$ be a graph with labeled vertices and edges. 
Let $f:G\rightarrow S^3$ be an embedding of $G$ 
into the oriented three sphere $S^3$. 
The embedding $f$ is called a {\it spatial graph}. 

Let $f_1$ and $f_2$ be spatial graphs. 
We say that $f_2$ is {\em obtained from $f_1$  
by a local move $(T_{1},T_{2})$} if 
there is an orientation preserving embedding 
$h:B^{3}\rightarrow S^{3}$ such that $f_{i}(G)\cap h(B^{3})=h(T_{i})$
for $i=1,2$ and $f_{1}(G- h(B^3))=
f_{2}(G-h(B^3))$ together with the labels of vertices and edges. 
Two spatial graphs $f_{1}$ and $f_{2}$ are {\it $C_{k}$-equivalent} 
if $f_{2}$ is obtained from $f_{1}$  by a finite sequence of $C_{k}$-moves 
and ambient isotopies. 
We note that the relation is an equivalence relation on 
spatial graphs. 
For a spatial graph $f$, let $[f]_k$ denote the $C_k$-equivalence 
class that contains $f$. 
It is known that $C_k$-equivalence implies 
$C_{k-1}$-equivalence \cite{Habiro1}, \cite{T-Y0}.

Let $l$ be a positive integer and $k_1,...,k_l$ positive integers. Suppose
that for each $P\subset\{1,...,l\}$, 
an spatial graph $f_P$ in $S^3$ is assigned. 
Suppose that there are mutually disjoint, orientation 
preserving embeddings 
$h_i:B^3\rightarrow S^3$ $(i=1,...,l)$ such that \\
(1) $f_{\emptyset}(G)-\bigcup_{i=1}^l h_i(B^3)=
f_{P}(G)-\bigcup_{i=1}^l h_i(B^3)$ together
with the labels 
for any subset $P\subset \{1,...,l\}$,\\
(2) $(h_i^{-1}(f_\emptyset(G)),h_i^{-1}(f_{\{1,...,l\}}(G)))$ is 
a $C_{k_i}$-move $(i=1,...,l)$, and\\
(3) $f_P(G)\cap h_i(B^3)=\left\{
\begin{array}{ll}
f_{\{1,...,l\}}(G)\cap h_i(B^3) & \mbox{if $i\in P$},\\
f_\emptyset(G)\cap h_i(B^3) & \mbox{otherwise}.
\end{array}
\right.$\\
Then we call the set $\{f_P|P\subset\{1,...,l\}\}$ 
a {\em singular spatial graph of type $(k_1,...,k_l)$}.

The {\em $\theta$-curve} 
is a graph $\theta$ with two vertices $v_1,v_2$ and 
three edges $e_1,e_2,e_3$  each of
which joins $v_1$ and $v_2$. 
When $G=\theta$, a spatial graph and singular spatial graph 
are called {\it spatial $\theta$-curve} and 
{\it singular spatial $\theta$-curve} respsctively.  
A spatial $\theta$-curve is {\em trivial} if its image is 
contained in a 2-sphere in $S^3$.
In the remainder of this sction, we consider only the 
case that a graph is the $\theta$-curve. 
 
Let $\Gamma$ be the set of all spatial $\theta$-curve types in
$S^3$ and ${\Bbb Z}{\Gamma}$
the free abelian group  generated by ${\Gamma}$. 
For a singular spatial $\theta$-curve 
${\bf f}=\{f_P|P\subset\{1,...,l\}\}$ of type
$(k_1,...,k_l)$, we define an element $\kappa({\bf f})$ of ${\Bbb Z}{\Gamma}$ by 
\[\kappa({\bf f})=\sum_{P\subset\{1,...,l\}}(-1)^{| P|}f_P.\]
Let ${\cal V}(k_1,...,k_l)$ be the subgroup of ${\Bbb Z}{\Gamma}$
generated by all $\kappa({\bf f})$ where $\bf f$ varies over all
singular spatial $\theta$-curves of type $(k_1,...,k_l)$.

For two spatial $\theta$-curves $f_1$ and $f_2$, 
remove small balls centred at $f_1(v_2)$ and $f_2(v_1)$ from $S^3$, 
then identify the boundaries so that the images of $i$-th edge are 
joined for each $i$. Then we obtain a new spatial $\theta$-curve. 
We call this embedding the {\em vertex connected sum} 
of $f_1$ and $f_2$, and denote by $f_1\# f_2$. 
The vertex connected sum is well-defined up to ambient 
isotopy \cite{Wolcott}. 
Let $f_1\# f_2$ be the vertex connected sum of two 
spatial $\theta$-curves $f_1$ and $f_2$. Then
$f_1\# f_2-f_1-f_2\in{\Bbb
Z}{\Gamma}$ is called a {\em composite relator}. 
Let ${\cal R}_\#$ be the subgroup of
${\Bbb Z}{\Gamma}$  generated by all composite relators.

Let $\iota:{\Gamma}\rightarrow{\Bbb Z}{\Gamma}$ be the natural 
inclusion map.
Let $\pi:{\Bbb Z}{\Gamma}\rightarrow
{\Bbb Z}{\Gamma}/{\cal V}(k_1,...,k_l)$ and $\lambda:
{\Bbb Z}{\Gamma}/{\cal V}(k_1,...,k_l)\rightarrow{\Bbb Z}{\Gamma}/({\cal
V}(k_1,...,k_l)+{\cal R}_\#)$ be the quotient homomorphisms. 
Then the composite maps
$\pi\circ\iota:
{\Gamma}\rightarrow {\Bbb Z}{\Gamma}/{\cal V}(k_1,...,k_l)$ and 
$\lambda\circ\pi\circ\iota:
{\Gamma}\rightarrow {\Bbb Z}{\Gamma}/({\cal V}(k_1,...,k_l)+{\cal R}_\#)$ 
are called the {\em universal Vassiliev invariant of type 
$(k_1,...,k_l)$} and {\em universal additive Vassiliev invariant of type 
$(k_1,...,k_l)$} respectively. We denote them by $v_{(k_1,...,k_l)}$ and
$w_{(k_1,...,k_l)}$
respectively.
In the case of knots, these are same invariants as defined 
by K. Taniyama and the author \cite{T-Y}.
Similarly, we can also define ${\cal V}(k_1,...,k_l)$ and 
the universal Vassiliev invariant 
$v_{(k_1,...,k_l)}$ for the embeddings of any graph. 
Since a $C_1$-move is a crossing change we see that a singular 
spatial graph of type
$(\underbrace{1,...,1}_{l})$ is essentially the same as a singular saptial 
graph with $l$ crossing vertices in the sense of T. Stanford \cite{Stanford}. 
Therefore we see that
$\displaystyle  v_{(1,...,1)}$ is the universal 
Vassiliev invariant of order $\leq l-1$. Note that $\displaystyle 
v_{(1,...,1)}(f_1)=\displaystyle 
v_{(1,...,1)}(f_2)$ if and only if $v(f_1)=v(f_2)$ for any Vassiliev
invariant $v$ of order $\leq l-1$. 
In the case of links, $v_{(2,...,2)}$ is the same as that 
defined in \cite{Mellor}, \cite{Stanford4}. 
In \cite{T-Y0} Taniyama and the author 
defined finite type invariants of order $(k;n)$ 
for the embeddings of a graph, 
which are essentially same as 
$\displaystyle v_{(\scriptsize\underbrace{n-1,...,n-1}_{k+1})}$.

By the arguments similar to that in Proofs of Theorems 1.1 and 1.2 
in \cite{T-Y} and that in Proof of Theorem 1.4 in \cite{T-Y}, 
we have the following two theorems.

\medskip
\noindent
{\bf Theorem 1.1.} {\em Let $k_1,...,k_l$ be positive integers 
and $k=k_1+\cdots+k_l$. Then the followings hold.

{\rm (1)} ${\cal V}(k)\subset {\cal V}(k_1,...,k_l)
\subset{\cal V}(\underbrace{1,...,1}_{k})$. 

{\rm (2)} ${\cal V}(k_1,...,k_l)+{\cal R}_\#=
{\cal V}(k)+{\cal R}_\#$. $\Box$}

\medskip
\noindent
{\bf Remark.} Theorem 1.1(1) holds for the spatial embeddings of any graph.

\medskip
\noindent
{\bf Theorem 1.2.} {\it The
$C_k$-equivalence classes of the spatial $\theta$-curves 
forms a group with the unit element $[f_0]_k$ 
under the vertex connected sum, where $f_0$ is a trivial 
$\theta$-curve. $\Box$}

\medskip
We denote  by $G_k$ this group. 
Let $\varphi:G_k\rightarrow{\Bbb Z}\Gamma/
{\cal V}(k)$ be a map induced by the 
inclusion $\iota:\Gamma\rightarrow{\Bbb Z}\Gamma$. 
By Theorem 1.2, $\varphi$ is a well-defined, 
epimorphism. (In fact, $\varphi$ is an isomorphism, see 
the remark after Corollary 1.4.) 
Since $\varphi([G_k,G_k])={\cal R}_\#$, by Theorem 1.1 (2), 
we have the following theorem.

\medskip
\noindent
{\bf Theorem 1.3.} {\it Let $k_1,...,k_l$ be positive integers 
and $k=k_1+\cdots+k_l$. Then $G_k/[G_k,G_k]$ is isomorphic to 
${\Bbb Z}{\Gamma}/({\cal V}(k_1,...,k_l)+{\cal R}_\#)$. $\Box$}

\medskip
Let $f_1$ and $f_2$ be spatial $\theta$-curves and $k=k_1+\cdots+k_l$. 
If $f_1$ and $f_2$ are $C_k$-equivalent, then by Theorem 1.1 (1), 
$f_1-f_2\in{\cal V}(k)\subset {\cal V}(k_1,...,k_l)$. 
Therefore we have $v_{(k_1,...,k_l)}(f_1)=v_{(k_1,...,k_l)}(f_2)$. 
On the other hand, if
$v_{(k_1,...,k_l)}(f_1)=v_{(k_1,...,k_l)}(f_2)$, then
$w_{(k_1,...,k_l)}(f_1)=w_{(k_1,...,k_l)}(f_2)$. Hence 
we have $f_1-f_2\in{\cal V}(k_1,...,k_l)+{\cal R}_\#$. 
If $G_k$ is abelian group, i.e., $[G_k,G_k]=\{id\}$, 
then  by Theorem 1.3, $f_1$ and $f_2$ are $C_k$-equivalent. 
So we have the following corollary. 

\medskip
\noindent
{\bf Corollary 1.4.} {\em Let $k_1,...,k_l$ be positive integers 
and $k=k_1+\cdots+k_l$. Let $f_1$ and $f_2$ be spatial $\theta$-curves.
If $G_k$ is an abelian group, 
then the following conditions are mutually equivalent.

{\rm (1)} $f_1$ and $f_2$ are $C_k$-equivalent,

{\rm (2)} $v_{(k_1,...,k_l)}(f_1)=v_{(k_1,...,k_l)}(f_2)$,

{\rm (3)} $w_{(k_1,...,k_l)}(f_1)=w_{(k_1,...,k_l)}(f_2)$.
$\Box$}

\medskip
\noindent
{\bf Remark.}
Let $f_1$ and $f_2$ be spatial graphs (not necessarily $\theta$-curve). 
If $f_1-f_2\in {\cal V}(k)$, then there are singular spatial graphs 
${\bf f}^i$'s of type $(k)$ and integers $x_i$'s such that 
$f_1-f_2= \sum x_i\kappa({\bf f}^i)$. By induction on $\sum|x_i|$, 
we see that $f_1$ and $f_2$ are $C_k$-equivalent. Since 
$f_1-f_2\in {\cal V}(k)$ if $f_1$ and $f_2$ are $C_k$-equivalent, 
we have the following: Two spatial graphs $f_1$ and $f_2$ are 
$C_k$-equivalent if and only if  $v_{(k)}(f_1)=v_{(k)}(f_2)$.

\medskip
\noindent
{\bf Theorem 1.5.} 
{\it Let $f_0$ be a trivial spatial $\theta$-curve.
Then the followings hold.

$(1)$ For each $f\in[f_0]_k$, $[f]_{k+2}$ belongs to the center of 
$G_{k+2}$.

$(2)$ If $2l\geq k$, then the set 
$H_k^l=\{[f]_k|f\in[f_0]_l\}$ is an abelian subgroup of $G_k$. }

\medskip
By \cite{tani} and \cite{M-T}, we have $[f_0]_1=[f_0]_2=\Gamma$. 
Hence, by Theorem 1.5(2), $H_4^2=G_4$ is abelian. 
By Corollary 1.4, we have

\medskip
\noindent
{\bf Corollary 1.6.} 
{\em Let $k_1,...,k_l$ $(l\leq 4)$be positive integers 
and $k=k_1+\cdots+k_l\leq 4$. Let $f_1$ and $f_2$ be spatial $\theta$-curves.
Then the following conditions are mutually equivalent.

{\rm (1)} $f_1$ and $f_2$ are $C_k$-equivalent,

{\rm (2)} $v_{(k_1,...,k_l)}(f_1)=v_{(k_1,...,k_l)}(f_2)$,

{\rm (3)} $w_{(k_1,...,k_l)}(f_1)=w_{(k_1,...,k_l)}(f_2)$.
$\Box$}

\medskip
\noindent
{\bf Remark.} As a special case of Corollary 1.6, we see that 
for $k\leq 4$ two spatial 
$\theta$-curves are $C_k$-equivalent if and only if the universal 
(additive) Vassiliev invariant of order $\leq k-1$ are equal. 
Meanwhile, a basis for the space of Vassiliev invariants of 
order $\leq 4$ is known \cite{Kan},\cite{Koi}.

\medskip
\noindent
{\bf Theorem 1.7.} {\em Let $f_0$ be a trivial spatial $\theta$-curve, and 
let $f_1$ and $f_2$ be in $[f_0]_k$. Then $f_1$ and $f_2$ are 
$C_{2k}$-equivalent  if and only if $v_{(k,k)}(f_1)=v_{(k,k)}(f_2)$.}

\medskip
As we saw before, $G_k$ is abelian for $k\geq 4$. However $G_k$ is not 
necessarily abelian. In fact, we have the following theorem. 

\medskip
\noindent
{\bf Theorem 1.8.} {\em The group $G_k$ is nonabelian for any 
$k\geq 12$.}


\medskip
\noindent
{\bf Remarks}
(1) If $G_k$ is abelian, then so is $G_{k'}$ for any $k'(<k)$. \\
(2) In the proof of Theorem 1.8, we see that there are two spatial 
$\theta$-curves $g$ and $h$ such that $g-h$ is in ${\cal R}_{\#}$ and 
is not in ${\cal V}(\underbrace{1,...,1}_{k})$ for $k\geq 12$. 
Hene $v_{(k_1,...,k_l)}(g)\neq v_{(k_1,...,k_l)}(h)$ for 
$k_1+\cdots +k_l\geq 12$ by Theorem 1.1(1), while 
$w_{(k_1,...,k_l)}(g)=w_{(k_1,...,k_l)}(h)$ for any $k_1,...,k_l$. 
In contrast, for any knots $K$ and $K'$, 
$v_{(k_1,...,k_l)}(K)=v_{(k_1,...,k_l)}(K')$ if and only if 
$w_{(k_1,...,k_l)}(K)=w_{(k_1,...,k_l)}(K')$ \cite{T-Y}.  \\
(3) Habiro showed the set of $C_k$-equivalence classes $S_k(m)$ of 
$m$-string links forms a group under the composition \cite{Habiro1}. 
By considering the complement of a regular neighbourhood of one of edges, 
we have that there is a surjection from the $2$-string links 
to the spatial $\theta$-curves. 
Since the surjection induces an epimorphism from $S_k(2)$ to $G_k$ 
and since there is an epimorphism from $S_k(m)\ (m>2)$ to $S_k(2)$, 
$S_k(m)$ is nonabelian for any $m\geq 2$ and $k\geq 12$.

\medskip
The following is still open. 

\medskip
\noindent
{\bf Problem.} Find the minimum number $k\ (5\leq k\leq 12)$ 
such that $G_k$ is nonabelian.

\bigskip\noindent
{\bf 2. Band description of spatial graphs}

\medskip
A {\it $C_{1}$-link model} is a pair $(\alpha,\beta)$ where 
$\alpha$ is a disjoint union of properly embedded arcs in $B^3$ and $\beta$ 
is a disjoint union of arcs on $\partial B^3$ with 
$\partial \alpha=\partial \beta$ as illustrated in Fig. 2.1.
Suppose that a $C_{k}$-link model $(\alpha, \beta)$ is defined 
where $\alpha$ is a disjoint union of $k+1$ properly embedded arcs in $B^3$ 
and $\beta$ is a disjoint union of $k+1$ arcs on $\partial B^3$ with 
$\partial \alpha=\partial \beta$ such that $\alpha \cup \beta$ is a disjoint 
union of $k+1$ circles. Let $\gamma$ be a component of $\alpha \cup \beta$ and 
$W$ a regular neighbourhood of $\gamma$ in $(B^3-(\alpha \cup
\beta))\cup\gamma$. Let $V$ be an
oriented solid  torus, $D$ a disk in $\partial V$, $\alpha_{0}$ properly
embedded arcs in $V$ 
and $\beta_{0}$ arcs on $D$ as illustrated in Fig. 2.2. Let $\psi
:V\rightarrow W$ be an orientation
preserving homeomorphism  such that $\psi (D)=W\cap \partial B^{3}$ and 
$\psi (\alpha_{0}\cup \beta_{0})$ bounds disjoint disks in $B^3$. 
Then we call the pair $((\alpha -\gamma )\cup \psi (\alpha_{0}), 
(\beta -\gamma )\cup\psi (\beta_{0}))$ a {\it $C_{k+1}$-link model.} 
A {\it link model} is a 
$C_{k}$-link model for some $k$. 
It is known that, for a $C_k$-link model $(\alpha,\beta)$, 
the local move 
$(\alpha,\hat{\beta})$ is equivalent to a $C_k$-move \cite{T-Y}, 
where $\hat{\beta}$ is a slight push in of $\beta$.

\begin{center} 
\begin{tabular}{cc}
\includegraphics[trim=0mm 0mm 0mm 0mm, width=.17\linewidth]
{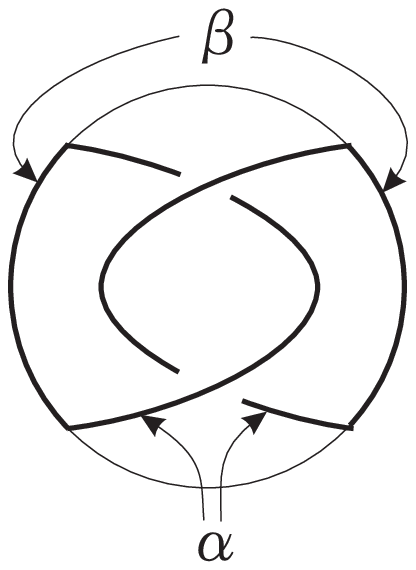} \hspace*{1in} &
\includegraphics[trim=0mm 0mm 0mm 0mm, width=.2\linewidth]
{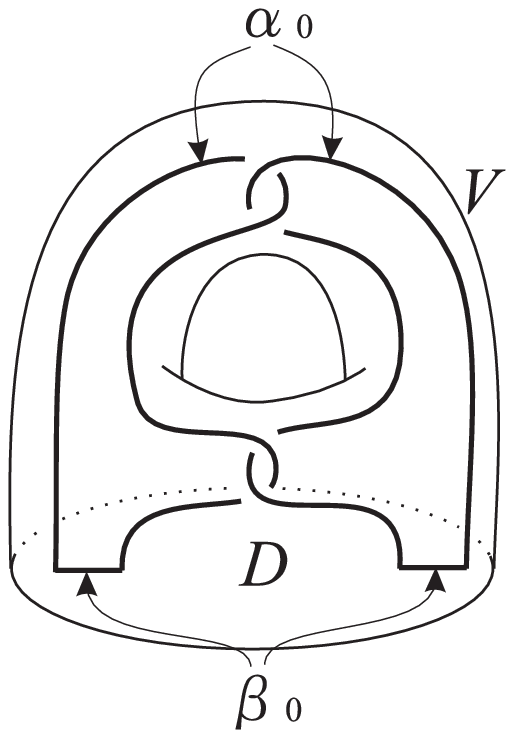}\\
Fig. 2.1 \hspace*{1in} &
Fig. 2.2
\end{tabular}
\end{center}

Let $f_1$ be a spatial $\theta$-curve, and let 
$(\alpha_{1},\beta_{1}),...,(\alpha_{l},\beta_{l})$ be link 
models. Let
$\psi_{i}:B^{3}\rightarrow
S^3$ $(i=1,...,l)$ be mutually disjoint, orientation preserving 
embeddings, and let 
$b_{1,1}, b_{1,2},..., b_{1,\rho(1)}, 
b_{2,1}, b_{2,2},..., b_{2,\rho(2)},
..., b_{l,1},b_{l,2},...,b_{l,\rho(l)}$
be mutually disjoint disks embedded in $S^3$. Suppose that 
they satisfy the following
conditions;\\
(1) $\psi_{i}(B^{3})\cap f_1(\theta)=\emptyset$  for each $i$,\\
(2) $b_{i,k}\cap f_1(\theta)=\partial b_{i,k}\cap 
(f_1(\theta)-f(v_1\cup v_2))$ 
is an arc for each $i,k$,\\
(3) $b_{i,k}\cap (\bigcup_{j=1}^{l} \psi_{j}(B^{3}))=
\partial b_{i,k}\cap \psi_{i}(B^{3})$ is a component of 
$\psi_{i}(\beta_{i})$ for each $i,k$,\\
(4) ($\bigcup_{k=1}^{\rho(i)}b_{i,k})\cap \psi_{i}(B^{3})
=\psi_{i}(\beta_{i})$ for each $i$.\\
Let $f_2$ be a spatial $\theta$-curve defined by
\[
f_2(\theta)=f_1(\theta)\cup (\bigcup_{i,k}\partial b_{i,k})\cup 
(\bigcup_{i=1}^{l}\psi_{i}(\alpha_{i})) - 
\bigcup_{i,k}{\rm int}(\partial b_{i,k}\cap f_1(\theta)) - 
\bigcup_{i=1}^{l}\psi_{i}({\rm int}\beta_{i}),
\]
where the labels of $f_2(\theta)$ 
coincides that of $f_1(\theta)$ on 
$f_1(\theta)-\bigcup_{i,k}b_{i,k}$. 
When $(\alpha_i,\beta_i)$ is a $C_{k}$-link model, we call 
$\psi_i(B^3)$ a {\it $C_k$-link ball}. 
We set 
${\cal B}_i=((\alpha_i,\beta_i),\psi_i,\{b_{i,1},...,b_{i,\rho(i)}\})$
and call ${\cal B}_i$ a {\it $C_k$-chord} when
$(\alpha_i,\beta_i)$ is a $C_k$-link model. 
We denote $f_2$  by
$\Omega(f_1;\{{\cal B}_1,...,{\cal B}_l\})$ and call it a {\it band
description} of $f_2$. We also say 
$f_2$ is a {\it band sum} of $f_1$ and 
link models $(\alpha_{1},\beta_{1}),...,(\alpha_{l},\beta_{l})$. 

By the arguments similar to that in Proof of Lemma 3.6 \cite{T-Y}, 
we have

\medskip
\noindent{\bf Lemma 2.1.} {\it 
Two spatial $\theta$-curves 
$f_1$ and $f_2$ are $C_k$-equivalence if and 
only if there are spatial $\theta$-curves $f'_1$ and $f'_2$ 
such that $f'_i$ is ambient isotopic to $f_i$ $(i=1,2)$ and 
$f'_2$ is a band sum of 
$f'_1$ and some $C_k$-link models. $\Box$}

\medskip
In the following lemma, the former assertion follows directly 
from Lemma 3.9 in \cite{T-Y} and the latter can be shown by the 
similar arguments as in proof of Lemma 3.9 in \cite{T-Y}. 

\medskip
\noindent{\bf Lemma 2.2.} {\it  A local move as illustrated 
in Fig. {\rm2.3} $($resp. Fig. {\rm 2.4}$)$ is realized  by a 
$C_{j+k}$-move $($resp. $C_{j+k+1}$-move$)$. $\Box$}

\begin{center}
\includegraphics[trim=0mm 0mm 0mm 0mm, width=.7\linewidth]
{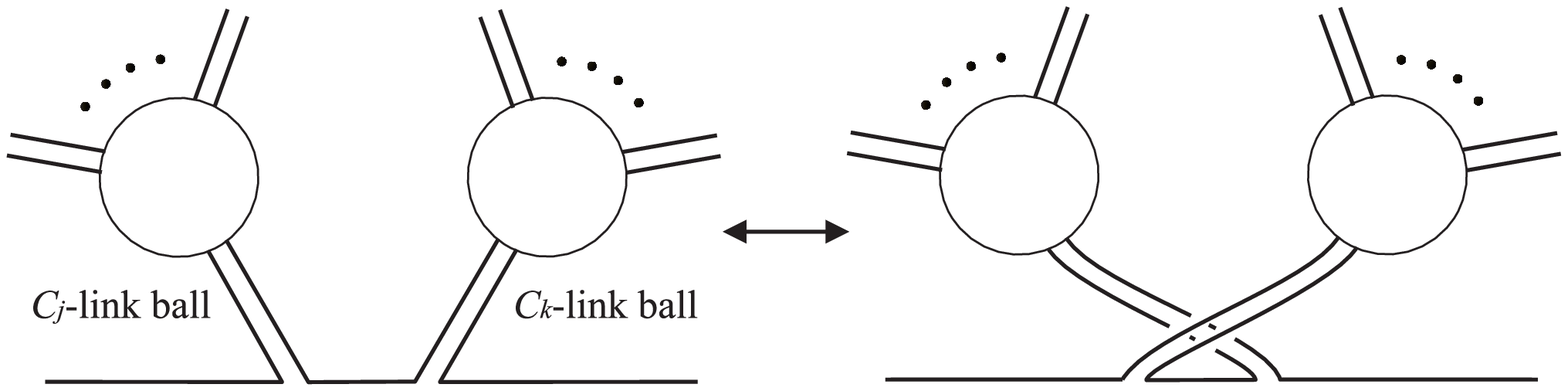}

Fig. 2.3
\end{center}

\begin{center}
\includegraphics[trim=0mm 0mm 0mm 0mm, width=.65\linewidth]
{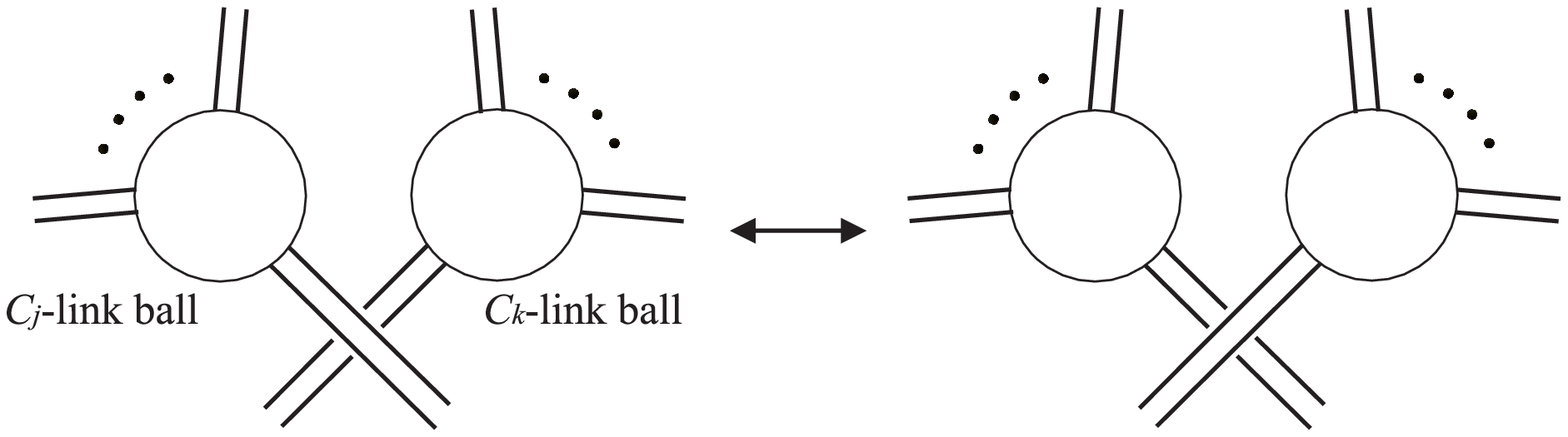}

Fig. 2.4
\end{center}

\medskip
\noindent
{\bf Proof of Theorem 1.5.} 
(1) Suppose $f\in[f_0]_k$, then by Lemma 2.1, 
we may assume that $f$ is a band sum of 
$f_0$ and some $C_k$-link models. 
Let $g$ be a spatial $\theta$-curve. 
Since $\Gamma=[f_0]_2$, we may suppose $g$  
is a band sum of $f_0$ and some $C_2$-link models. 
It is not hard to see that $f\# g$ and $g\# f$ are 
transposed each other by the moves as in Figs. 2.3 and 
2.4, where we consider the case $j=2$, and ambient isotopies. 
Thus by Lemma 2.2 we have $f\# g$ is $C_{k+2}$-equivalent to 
$g\# f$. Hence we have
$[f]_{k+2}[g]_{k+2}=[f\# g]_{k+2}=[g\# f]_{k+2}=[g]_{k+2}[f]_{k+2}$. 

(2) Suppose that both $[f_1]_k$ and $[f_2]_k$ belong to $H_k^l$. 
Then we note that $[f_1]_k[f_2]_k=[f_1\# f_2]_k\in H_k^l$. 
If $[f]_k$ belongs to $H_k^l$, then by Lemma 2.1, we may assume that 
$f=\Omega(f_0;\{{\cal B}_1,...,{\cal B}_n\})$ and 
$f_0=\Omega(f;\{{\cal B}'_1,...,{\cal B}'_m\})$ for some 
$C_l$-chords ${\cal B}_1,...,{\cal B}_n,{\cal B}'_1,...,{\cal B}'_m$. 
By using Sublemma 3.5 in \cite{T-Y} repeartedly, there are $C_l$-chords 
${\cal B}''_1,...,{\cal B}''_m$ such that 
$f_0$ is ambient isotopic to $\Omega(f_0;\{{\cal B}_1,...,{\cal B}_n
{\cal B}''_1,...,{\cal B}''_m\})$. 
By Lemma 2.2, we can deform 
$\Omega(f_0;\{{\cal B}_1,...,{\cal B}_n
{\cal B}''_1,...,{\cal B}''_m\})$ into 
$\Omega(f_0;\{{\cal B}_1,...,{\cal B}_n\})\#
\Omega(f_0;\{{\cal B}''_1,...,{\cal B}''_m\})$
by $C_{2l}$-moves and ambient isotopies, i.e., $f_0$ and $f\# 
\Omega(f_0;\{{\cal B}''_1,...,{\cal B}''_m\})$ are 
$C_{2l}$-equivalent. Since $2l\geq k$, 
$f_0$ and $f\# \Omega(f_0;\{{\cal B}''_1,...,{\cal B}''_m\})$ are 
$C_{k}$-equivalent. This implies that 
$[\Omega(f_0;\{{\cal B}''_1,...,{\cal B}''_m\})]_k=[f]_k^{-1}
\in H_k^l$. Therefore $H_k^l$ is a subgroup of $G_k$. 
By the arguments similar to that in (1), we see that 
$H_k^l$ is abelian. 
$\Box$

\medskip
Let $v$ be an invariant of the embeddings of a graph that takes 
values in an abelian group. We call $v$ a {\em Vassiliev invariant of 
type $(k_1,...,k_l)$} if, for any singular spatial graph 
$\{f_P|P\subset\{1,...,l\}\}$ of type $(k_1,...,k_l)$, 
\[\sum_{P\subset \{1,...,l\}}(-1)^{|P|}v(f_P)=0.\]

\medskip
\noindent
{\bf Proof of Theorem 1.7.} The \lq only if'  part follows from 
Theorem 1.1(1). We shall show \lq if'  part. Let 
$\varphi:\Gamma\rightarrow H_{2k}^k$ be a map defined as follows: 
\[\varphi(f)=
\left\{\begin{array}{ll}
[f]_{2k}&\mbox{if }f\in[f_0]_k,\\
0&\mbox{otherwise}.
\end{array}\right.\]
Clearly, $\varphi$ is an invarinat. 
Now we will show that $\varphi$ is a Vassiliev invariant 
of type $(k,k)$. 
Let $\{h_P|P\subset\{1,2\}\}$ be a singular spatial $\theta$-curve 
of type $(k,k)$. Since $\{h_P|P\subset\{1,2\}\}\subset[h_{\emptyset}]_k$, 
we have $\varphi(h_P)=0$ if $[h_{\emptyset}]_k\neq[f_0]_k$. 
So we may suppose that $[h_{\emptyset}]_k=[f_0]_k$. 
Then we have 
\[\sum_{P\subset \{1,2\}}(-1)^{|P|}\varphi(h_P)=
[h_{\emptyset}]_{2k}-[h_{\{1\}}]_{2k}
-[h_{\{2\}}]_{2k}+[h_{\{1,2\}}]_{2k}.\]
Since $[h_{\emptyset}]_k=[f_0]_k$, by Lemma 2.1, we may assume 
that $h_{\emptyset}=\Omega(f_0;\{{\cal B}_1,...,{\cal B}_n\})$ for 
some $C_k$-chords ${\cal B}_1,...,{\cal B}_n$. 
By Sublemma 3.1 in \cite{T-Y} (or Lemma 3.7 in the next section), 
there are $C_k$-chords ${\cal B}'_1,{\cal B}'_2$ such that 
$\Omega(h_{\emptyset};\cup_{i\in P}\{{\cal B}'_i\})$ is ambient isotopic 
to $h_P$. By the arguments similar to that in the proof of Sublemma 3.5 in 
\cite{T-Y}, we see that there are $C_k$-chords ${\cal B}''_1,{\cal B}''_2$ 
such that 
$\Omega(f_0;\{{\cal B}_1,...,{\cal B}_n\}\cup(\cup_{i\in P}\{{\cal B}''_i\}))$ 
is ambient isotopic to $h_P$. 
Since ${\cal B}_1,...,{\cal B}_n, {\cal B}''_1,{\cal B}''_2$ are 
$C_k$-chords, by Lemma 2.2, we have 
\[
[\Omega(f_0;\{{\cal B}_1,...,{\cal B}_n\}\cup
\{{\cal B}''_1,{\cal B}''_2\})]_{2k}=
[\Omega(f_0;\{{\cal B}_1,...,{\cal B}_n\})\#
\Omega(f_0;\{{\cal B}''_1\})\#\Omega(f_0;\{{\cal B}''_2\})]_{2k},\]  
and 
\[[\Omega(f_0;\{{\cal B}_1,...,{\cal B}_n\}\cup \{{\cal B}''_i\})]_{2k}
=[\Omega(f_0;\{{\cal B}_1,...,{\cal B}_n\})\# 
\Omega(f_0;\{{\cal B}''_i\})]_{2k}\ (i=1,2).\] 
Hence we have 
\[\begin{array}{l}
\sum_{P\subset \{1,2\}}(-1)^{|P|}\varphi(h_P)\\
\hspace*{2.5cm}=
[f_{\emptyset}]_{2k}-[f_{\emptyset}\#\Omega(f_0;\{{\cal B}''_1\})]_{2k}
-[f_{\emptyset}\#\Omega(f_0;\{{\cal B}''_2\})]_{2k}\\
\hspace*{7cm}+[f_{\emptyset}\#\Omega(f_0;\{{\cal B}''_1\})
\#\Omega(f_0;\{{\cal B}''_2\})]_{2k}=0\in H_{2k}^k.
\end{array}\]
Therefore $\varphi$ is a Vassiliev invariant of type $(k,k)$. 
This and the assumption $v_{(k,k)}(f_1)=v_{(k,k)}(f_2)$ imply  
$\varphi(f_1)=\varphi(f_2)$. By the definition of $\varphi$, we have 
$[f_1]_{2k}=[f_2]_{2k}$.  $\Box$

\bigskip\noindent
{\bf 3. Disk/band surfaces and Vassiliev invariants of spatial graphs}

\medskip
A graph $G$ is {\it trivalent} if the valence of any vertex of $G$ is 
equal to $3$. 
A graph $G$ is {\it planar} if there exists an embedding 
$f_0:G\rightarrow  {\Bbb R}^2$.
A connected, planar graph $G$ is said to be {\it prime} if, for any 
embedding $f_0:G\rightarrow  {\Bbb R}^2$, there exist no 
simple closed 
curves $C$ in ${\Bbb R}^2$ satisfying either the following (1) or (2) 
(cf. \cite{suzuki}, \cite{KSWZ}), where $A$, $B$ are the two components of 
${\Bbb R}^2-C$.\\
(1) $C$ meets $f_0(G)$ in a single point such that both 
$A\cap f_0(G)$ and $B\cap f_0(G)$ are non-empty.\\
(2) $C$ meets $f_0(G)$ in two points such that both 
$A\cap f_0(G)$, $B\cap f_0(G)$ are neither empty nor single open arcs.

For any connected, planar graph $G$, we fix a planar embedding 
$f_0:G\rightarrow  {\Bbb R}^2$ arbitrarily.
The image $f_0(G)$ has complementary domains $D_1,D_2,...,D_n$ 
that are bounded and one unbounded $D_0$.
The preimage $c_i=f^{-1}_0(\partial D_i)$ is a 1-complex which 
can be viewed as a 1-cycle in $H_1(G;{\Bbb Z})$.
We call $c_i$ $(i\neq 0)$, $c_0$ respectively a {\em boundary cycle} 
and the {\em outermost cycle} in $G$ with respect to $f_0$.

For a spatial embedding $f:G\rightarrow  S^3$ of a graph $G$, a 
{\em disk/band surface} $S$ of $f(G)$ is a compact, 
orientable surface in $ S^3$ such that $f(G)$ is a deformation 
retract of $S$ contained in ${\rm int}S$ \cite{KSWZ}. 

In \cite{SSY}, T. Soma, H. Sugai and the author showed the 
following theorem.

\medskip
\noindent
{\bf Theorem 3.1.} (\cite[Theorem 1]{SSY}) 
{\em Suppose that $G$ is a connected, planar, 
prime and trivalent graph, 
and $f_0:G\rightarrow  {\Bbb R}^2$ is an embedding. 
Then, for any embedding $f:G\rightarrow   S^3$, 
there exists the unique disk/band surface $S$ of $f(G)$ 
up to ambient isotopy of which the Seifert pairings satisfying the 
following equation.
\[\langle f(c_i),f(c_j)\rangle_S
=\left\{ 
\begin{array}{ll}
-{\mathrm lk}(f(c_i),f(c_0))-\!\! \displaystyle{\sum_{c_i\cap c_k=\emptyset}}
\!\! {\mathrm lk}(f(c_i),f(c_k)) & \mbox{\sl if $i=j$ and $c_i\cap c_0=
\emptyset$}\\
0 & \mbox{\sl if $i\neq j$ and $c_i\cap c_j\neq\emptyset$}\\[2mm]
0 & \mbox{\sl if $i= j$ and $c_i\cap c_0\neq\emptyset$}\\[2mm]
{\mathrm lk}(f(c_i),f(c_j)) & \mbox{\sl if  $c_i\cap c_j=\emptyset$},
\end{array}
\right. \]
where $c_i,\ c_j,\ c_k$ are boundary cycles and $c_0$ is 
the outermost cycle with respect to $f_0$.  $\Box$}

\medskip
\noindent
We call the disk/band surface above the {\em canonical disk/band 
surface} for $f$. 
Note that the Seifert linking form of the canonical disk/band surface 
depends only on the linking numbers of pairs of disjoint cycles. 
If $G$ is the $\theta$-curve or the complete graph 
with $4$ vertices, then the cannonical disk/band surface is same as 
the {\em disk/band surface with zero Seifert linking form} 
that is defined in \cite{KSWZ}.
By the proof of Theorem 1 in \cite{SSY}, we note that 
the canonical disk/band surface is given as the image of an embedding 
of the regular neighborhood $S_0$ of $f_0(G)$ in ${\Bbb R}^2$. 
Thus by fixing orientation and label of $\partial S_0$, we 
have an ordered, oriented link as the image of 
an embedding of $\partial S_0$. 
From now on, we always assume that, for each graph $G$, 
$\partial S_0$ has fixed orientation and label.

Let $(T_{1},T_{2})$ be a 
local move, $t_{11},...,t_{1n}$ the components of $T_1$ and 
$t_{21},...,t_{2n}$ the components of $T_2$ 
with $\partial t_{1i}=\partial t_{2i}$ $(i=1,...,n)$. 
Let $N_{1i}$ and $N_{2i}$ be regular 
neighbourhoods of $t_{1i}$ and $t_{2i}$ in $B^3$ respectively such
that $N_{1i}\cap \partial  B^{3}=N_{2i}\cap \partial B^{3}$ $(i=1,...,n)$ 
and $N_{1i}\cap N_{1j}=N_{2i}\cap N_{2j}=\emptyset$ $(1\leq i<j\leq n)$. 
Let $\alpha_i$ $(i=1,...,n)$ be disjoint union of properly 
embedded $l_i$ arcs in $B^{2}\times [0,1]$ as illustrated in Fig. 3.1. 
Let $\psi_{ji}:B^{2}\times [0,1]\rightarrow N_{ji}$ be a homeomorphism 
with $\psi_{ji}(B^{2}\times \{ 0,1\} )=N_{ji}\cap \partial B^{3}$ for 
$j=1,2,\ i=1,...,n$. 
Suppose that $\psi_{1i}(\partial \alpha_i )=\psi_{2i}(\partial \alpha_i )$ and 
$\psi_{1i}(\alpha_i )$ and $\psi_{2i}(\alpha_i )$ are ambient isotopic in $B^{3}$ 
relative to $\partial B^3$. Then we say that a local move 
$(\bigcup_{i=1}^n \psi_{1i}(\alpha_i ), \bigcup_{i=1}^n \psi_{2i}(\alpha_i))$ 
is a {\it parallel of $(T_{1},T_{2})$ with weight} $(l_1,...,l_n)$.

\begin{center}
\includegraphics[trim=0mm 0mm 0mm 0mm, width=.25\linewidth]
{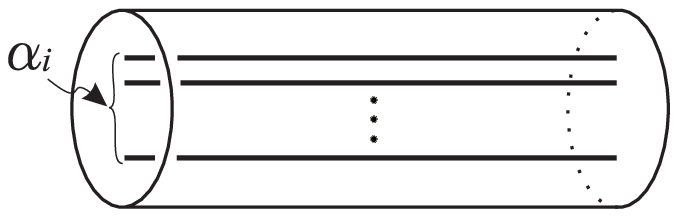}

Fig. 3.1
\end{center}

\medskip
\noindent
{\bf Proposition 3.2.} {\em Let $G$ be a connected, planar, prime 
and trivalent graph. Let $f_i:G\rightarrow S^3$ $(i=1,2)$ be 
embeddings and $S_i$ the canonical disk/band surface for $f_i$ $(i=1,2)$. 
If $f_1$ and $f_2$ are $C_k$-equivalent, then 
$\partial S_1$ and $\partial S_2$ are $C_k$-equivalent. }

\medskip
Let $T_1$ and $T_2$ be tangles. We say that $T_2$ is 
{\em obtained from $T_1$ by a local move} $(U_1,U_2)$ 
if there is an orientation preserving embedding 
$h:B^3\rightarrow{\mathrm int}B^3$ such that 
$S_i\cap h(B^3)=h(U_i)$ for $i=1,2$ and 
$T_1-h(B^3)=T_2-h(B^3)$. 
Two tangles $T_1$ and $T_2$ are {\it $C_k$-equivalent} 
if $T_2$ is obtained from $T_1$ by a finite sequence of 
$C_k$-moves and ambient isotopies relative $\partial B^3$. 

\medskip
\noindent
{\bf Lemma 3.3.} (cf. \cite[Claim on p. 26]{Habiro1}) 
{\em Let $(T_1,T_2)$ be a parallel of a $C_k$-move. 
Then $T_1$ and $T_2$ are $C_k$-equivalent. }

\medskip
\noindent
{\bf Proof.} 
Let $(U_1,U_2)$ be a $C_k$-move, $u_{11},...,u_{1k+1}$ and 
$u_{21},...,u_{2k+1}$ the components of $U_1$ and $U_2$ 
respectively. Suppose that $(T_1,T_2)$ is a parallel of 
$(U_1,U_2)$ with weight $(l_1,...,l_{k+1})$. 
We give a proof by induction on $l_1\times\cdots\times l_{k+1}$. 
In the case that $l_1\times\cdots\times l_{k+1}=1$, it obviously holds. 
Suppose $l_1\times\cdots\times l_{k+1}\geq 2$. We may suppose $l_1\geq 2$. 
By Lemma 2.1 in \cite{T-Y}, we may assume that the $C_k$-move 
$(U_1,U_2)$ is as illustrated in Fig. 3.2, i.e., 
the arcs except for $u_{i1}$ are contained in the shaded part in 
$U_i$ $(i=1,2)$.  
It is not hard to see that $T_2$ is obtained from $T_1$ by 
$l_1$ local moves that are paralells of $(U_1,U_2)$ with weight 
$(1,l_2,...,l_{k+1})$. This completes the proof. $\Box$

\begin{center}
\includegraphics[trim=0mm 0mm 0mm 0mm, width=.35\linewidth]
{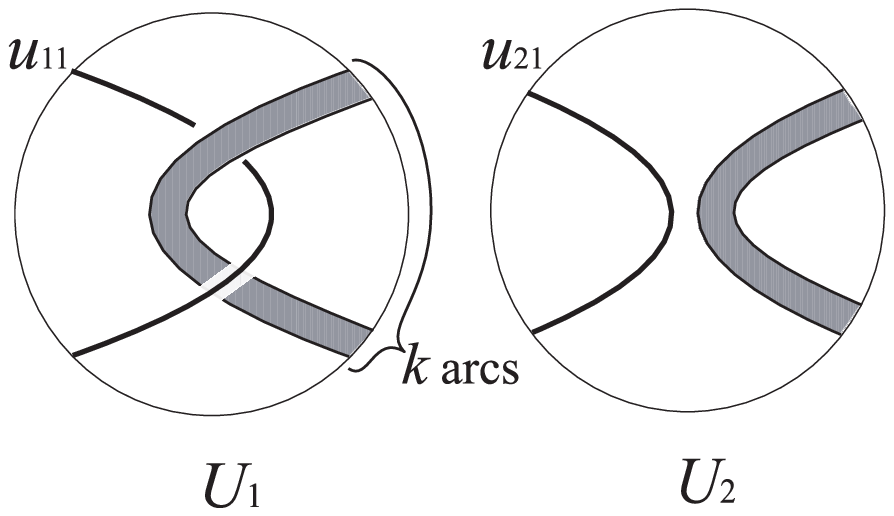}

Fig. 3.2
\end{center}

\medskip
\noindent
{\bf Proof of Proposition 3.2.} 
In the case $k=1$, it clearly holds. 
We consider the case $k\geq 2$. 
It is sufficient to consider the case that $f_2$ is obtained 
from $f_1$ by a single $C_k$-move. 
Suppose that $f_2$ is obtained from $f_1$ by a $C_k$-move. 
Then there is an embedding $h:B^3\rightarrow S^3$ such that 
$(h^{-1}(f_1(\theta)),h^{-1}(f_2(\theta)))$ is a $C_k$-move and 
$f_1(\theta)- h(B^3)=f_2(\theta)- h(B^3)$ together with the labels. 
We may suppose that 
$(h^{-1}(\partial S_1),h^{-1}(\partial S_2))$ is a parallel 
of the $C_k$-move $(h^{-1}(f_1(\theta)),h^{-1}(f_2(\theta)))$ 
with weight $(2,...,2)$. 
Then we have a new disk/band surface $S'_2$ for $f_2$ from $S_1$ 
by the local move 
$(h^{-1}(\partial S_1),h^{-1}(\partial S_2))$. 
Since $C_k$-move ($n\geq 2$) does not change the linking number, 
by Theorem 3.1, we have 
$\langle f_1(c_i),f_1(c_j)\rangle_{S_1}=
\langle f_2(c_i),f_2(c_j)\rangle_{S_2}$, and since 
the local move 
$(h^{-1}(\partial S_1),h^{-1}(\partial S_2))$ 
does not change the Seifert linking form of a disk/band surface, 
we have $\langle f_1(c_i),f_1(c_j)\rangle_{S_1}=
\langle f_2(c_i),f_2(c_j)\rangle_{S'_2}$. 
So we have $\langle f_2(c_i),f_2(c_j)\rangle_{S_2}=
\langle f_2(c_i),f_2(c_j)\rangle_{S'_2}$. 
By Theorem 3.1, $S_2$ is ambient isotopic to $S'_2$. 
Thus $\partial S_2$ is obtained from $\partial S_1$ by a parallel of 
a $C_k$-move. Lemma 3.3 completes the proof. $\Box$ 

\medskip
Let $G$ be a connected, planar, prime and trivalent graph and 
$\Gamma(G)$ the set of spatial graph types. 
Let $v$ be an invariant of ordered, oriented links that takes 
values in an abelian group $A$. 
Then we define a map $s:\Gamma(G)\rightarrow A$ as 
$s(f)=v(\partial S)$, where $S$ is the canonical disk/band 
surface for $f$. By Theorem 3.1, $s$ is an invariant of $f$. 
We call $s$ the {\em invariant induced from} $v$.  

\medskip
\noindent
{\bf Theorem 3.4.} 
{\em Let $G$ be a connected, planar, prime and trivalent graph and 
$\Gamma(G)$ the set of spatial graph types. 
Let $v$ be a Vassiliev invariant of 
type $(k_1,...,k_l)$ for ordered, oriented links. 
Then the invariant for $\Gamma(G)$ induced from $v$ is 
a Vassiliev invariant of type $(k_1,...,k_l)$.  }

\medskip
\noindent
In Theorem 3.4, the case of that a graph is the $\theta$-curve and 
$k_1=\cdots =k_l=1$ is given by Stanford \cite{Stanford0}.

By the arguments similar to that in the proof of Lemma 1.4 in 
\cite{Stanford0}, we have

\medskip
\noindent
{\bf Lemma 3.5.}
{\em Let $v$ be a Vassiliev invariant of type $(k_1,...,k_l)$ for 
ordered, oriented links and $s$ the invariant induced from $v$. 
Let $\{f_P|P\subset\{1,...,l\}\}$ be a singular spatial graph 
of type $(k_1,...,k_l)$. Let $L_P$ be the ordered, oriented link 
that is the boundary of the canonical disk/band surface for $f_P$ 
$(P\subset\{1,...,l\})$. 
Suppose there are mutually disjoint embeddings   
$h_{ij}:B^3\rightarrow S^3$ $(i=1,...,l,\ j=1,...,n_i)$ such that \\
$(1)$ $L_{\emptyset}-\bigcup_{i,j} h_{ij}(B^3)=
L_{P}-\bigcup_{i,j} h_{ij}(B^3)$ together
with orientations and labels of the components 
for any subset $P\subset \{1,...,l\}$,\\
$(2)$ $(h_{ij}^{-1}(L_\emptyset),
h_{ij}^{-1}(L_{\{1,...,l\}}))$ is 
a $C_{k_i}$-move $(i=1,...,l,\ j=1,...,n_i)$, and\\
$(3)$ $L_P\cap h_{ij}(B^3)=\left\{
\begin{array}{ll}
L_{\{1,...,l\}}\cap h_{ij}(B^3) & \mbox{if $i\in P$},\\
L_\emptyset\cap h_{ij}(B^3) & \mbox{otherwise}.
\end{array}
\right.$\\
Then we have 
\[\sum_{P\subset \{1,...,l\}}(-1)^{|P|}s(f_P)
=\sum_{P\subset \{1,...,l\}}(-1)^{|P|}v(L_P)=0. \Box \]
}

The following lemma follows directly from 
the proof of Theorem 1 in \cite{SSY}. 

\medskip
\noindent
{\bf Lemma 3.6.}
{\em Let $G$ be a connected, planar, prime and trivalent graph,  
$f_0:G\rightarrow{\Bbb R}^2$ an embedding, and 
$S_0$ the regular neighborhood of $f_0(G)$ in ${\Bbb R}^2$. 
Let $S$ be a disk/band surface for an embedding $f$. Suppose that 
$S$ is the image of an embedding of $S_0$ that is an extension 
of $f$. Then the canonical disk/band surface for $f$ is obtained from $S$ 
by a finite sequence of the moves as illustrated in 
Fig. $3.3$. $\Box$}

\begin{center}
\includegraphics[trim=0mm 0mm 0mm 0mm, width=.3\linewidth]
{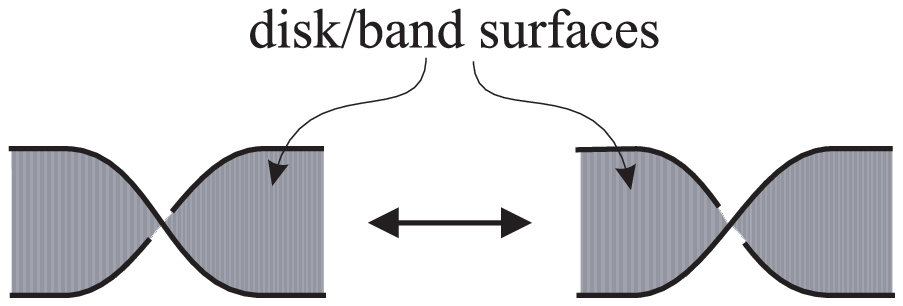}

Fig. 3.3
\end{center}

In the definition of band sum in Section 2, 
by replacing $f_i(\theta)$ with $T_i$ $(i=1,2)$, 
we can define that $T_2$ is a {\it band sum} of $T_1$ and link models 
$(\alpha_{1},\beta_{1}),...,(\alpha_{l},\beta_{l})$. 
By the arguments similar to that in Proof of Lemma 3.6 \cite{T-Y}, 
we have the following lemma. 

\medskip
\noindent
{\bf Lemma 3.7.} {\em 
Two tangles $T_1$ and $T_2$ are $C_k$-equivalent if and 
only if there are tangless $T'_1$ and $T'_2$ 
such that $T'_i$ is ambient isotopic to $T_i$ $(i=1,2)$ 
relative $\partial B^3$ and 
$T'_2$ is a band sum of 
$T'_1$ and some $C_k$-link models. $\Box$
}

\medskip
\noindent
{\bf Lemma 3.8.} {\em
Let $T_1$ and $T_2$ be tangles. If $T_2$ is  obtained from 
$T_1$ by a parallel of a $C_k$-move, then there are tangles 
$T'_1,T'_2$ and mutually disjoint, orientation preserving 
embeddings $h_i:B^3\rightarrow{\mathrm int}B^3$ $(i=1,...,n)$ 
such that\\
$(1)$ $T'_j$ is ambient isotopic to $T_j$ $(j=1,2)$ 
relative $\partial B^3$,\\
$(2)$ $T'_1-\bigcup_i h_i(B^3)=T'_2-\bigcup_i h_i(B^3)$, and\\
$(3)$ $(h_i^{-1}(T'_1),h_i^{-1}(T'_2))$ is a $C_k$-move $(i=1,...,n)$.}

\medskip
\noindent
{\bf Proof.} 
By Lemma 3.3, $T_1$ and $T_2$ are $C_k$-equivalent. 
By Lemma 3.7, 
there are tangles $T'_1$ and $T'_2$ such that 
$T'_j$ is ambient isotopic to $T_j$ $(j=1,2)$ relative $\partial B^3$ 
and that $T'_2$ is a band sum of $T'_1$ and some $C_k$-link models 
$(\alpha_1,\beta_1),...,(\alpha_n,\beta_n)$. 
Since $(\alpha_i,\hat{\beta}_i)$ $(i=1,...,n)$ are $C_k$-moves, 
we have the conclusion. $\Box$

\medskip
\noindent
{\bf Proof of Theorem 3.4.}
Since $\{f_P|P\subset \{1,...,l\}\}$ is a singular spatial graph 
of type $(k_1,...,k_l)$, by the definition, 
there are mutually disjoint, orientation 
preserving embeddings 
$h_i:B^3\rightarrow S^3$ $(i=1,...,l)$ such that \\
(1) $f_{\emptyset}(G)-\bigcup_{i} h_i(B^3)=
f_{P}(G)-\bigcup_{i} h_i(B^3)$ together
with the labels for any subset $P\subset \{1,...,l\}$,\\
(2) $(h_i^{-1}(f_\emptyset(G)),h_i^{-1}(f_{\{1,...,l\}}(G)))$ is 
a $C_{k_i}$-move $(i=1,...,l)$, and\\
(3) $f_P(G)\cap h_i(B^3)=\left\{
\begin{array}{ll}
f_{\{1,...,l\}}(G)\cap h_i(B^3) & \mbox{if $i\in P$},\\
f_\emptyset(G)\cap h_i(B^3) & \mbox{otherwise}.
\end{array}
\right.$\\
Let $S_{\emptyset}$ be the canonical disk/band surface for 
$f_{\emptyset}$.  
By considering the intersections $S_{\emptyset}\cap 
h_{i}(B^3)$ $(i=1,...,l)$, we find 
disk/band surfaces $S_P$ for $f_P$ $(P\subset \{1,...,l\})$ 
such that\\
(1) $S_{\emptyset}-\bigcup_{i} h_{i}(B^3)=
S_{P}-\bigcup_{i} h_{i}(B^3)$\\
(2) $(h_{i}^{-1}(\partial S_\emptyset),
h_{i}^{-1}(\partial S_{\{1,...,l\}}))$ is 
a parallel of a $C_{k_i}$-move 
$(i=1,...,l)$, and\\
(3) $S_P\cap h_{i}(B^3)=\left\{
\begin{array}{ll}
S_{\{1,...,l\}}\cap h_{i}(B^3) & \mbox{if $i\in P$},\\
S_\emptyset\cap h_{i}(B^3) & \mbox{otherwise}.
\end{array}
\right.$\\
By the proof of Proposition 3.2, if $k_i\geq 2$ for any $i\in P$, 
then $S_P$ is the canonical disk/band surface.   
Set $h_i=h_{i1}$ $(i=1,...,l)$. 
By Lemma 3.6, there are mutually disjoint, orientation 
preserving embeddings $h_{ij}:B^3\rightarrow S^3$ 
$(i=1,...,l,\ j=1,...,n_i)$, where $n_i=1$ if $k_i\geq 2$, 
and the canonical disk/band surfaces 
$S'_P$ for $f_P$ such that\\
(1) $S'_{\emptyset}-\bigcup_{i,j} h_{ij}(B^3)=
S'_P-\bigcup_{i,j} h_{ij}(B^3)$,\\
(2) $(h_{ij}^{-1}(\partial S'_{\emptyset}), 
h_{ij}^{-1}(\partial S'_{\{1,...,l\}}))$ is a paralell of 
$C_{k_i}$-move $(i=1,...,l,\ j=1,...,n_i)$, and\\
(3) $S'_P\cap h_{ij}(B^3)=\left\{
\begin{array}{ll}
S'_{\{1,...,l\}}\cap h_{ij}(B^3) & \mbox{if $i\in P$},\\
S'_\emptyset\cap h_{ij}(B^3) & \mbox{otherwise}.
\end{array}
\right.$\\
By combining this, Lemmas 3.8 and 3.5, we have the conclusion. 
$\Box$

\medskip
Let $G$ be a connected, planar, prime and trivalent graph 
and $E(G)=\{e_1,...,e_n\}$ the set of edges of $G$. 
Let $S_f$ be the cannonical disk/band surface for a spatial embedding 
$f$ of $G$, and let $S_f(x_1,...,x_n;y_1,...,y_n)$ 
$(x_i\in{\Bbb Z},y_j\in\{-1,0,1\})$ be a surface obtained from
$S_f$ as illustrated in Fig. 3.4. We note that 
$S_f(x_1,...,x_n;y_1,...,y_n)$ depends only on $S_f$ and the 
integers $x_1,...,x_n,y_1,...,y_n$. This means 
$S_f(x_1,...,x_n;y_1,...,y_n)$ is the unique surface for 
$f$. Let $v$ be an invariant of ordered, oriented links that takes 
values in an abelian group $A$. 
Then we can define an invariant 
$s_{(x_1,...,x_n;y_1,...,y_n)}:\Gamma(G)\rightarrow A$ as 
$s_{(x_1,...,x_n;y_1,...,y_n)}(f)=v(\partial S_f(x_1,...,x_n;y_1,...,y_n))$. 
We call $s_{(x_1,...,x_n;y_1,...,y_n)}$ the {\em invariant induced from} $v$ 
with respect to $x_1,...,x_n,y_1,...,y_n$.

\begin{center}
\includegraphics[trim=0mm 0mm 0mm 0mm, width=.7\linewidth]
{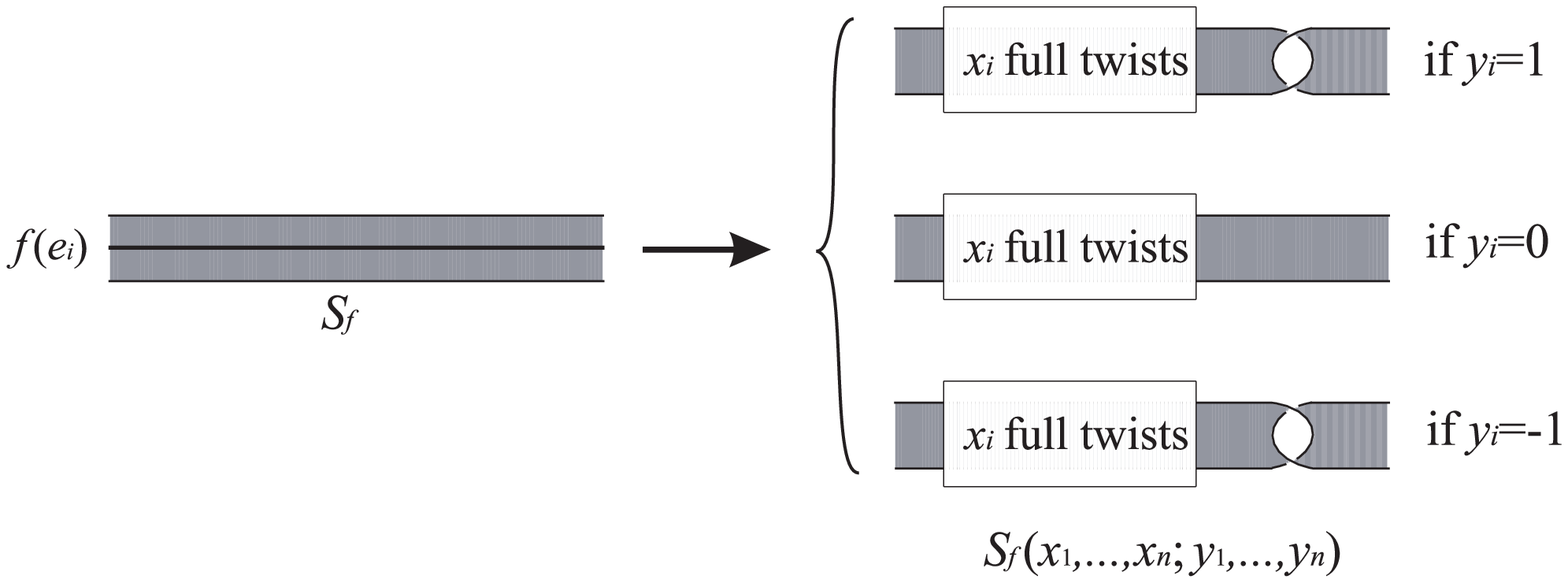}

Fig. 3.4
\end{center}

\noindent
By the arguments similar to that in the proofs of Proposition 3.2 and 
Theorem 3.4, we have  the following theorem.

\medskip
\noindent
{\bf Theorem 3.9} {\em Let $G$ be a connected, planar, prime 
and trivalent graph and $E(G)=\{e_1,...,e_n\}$ the set of 
edges of $G$. Then the followings hold.

{\rm (1)} Let $f_1$ and $f_2$ be spatial graphs and 
$S_{i}(x_1,...,x_n;y_1,...,y_n)$ the surface 
obtained from the canonical disk/band surface for $f_i$ $(i=1,2)$. 
If $f_1$ and $f_2$ are $C_k$-equivalent, then 
$\partial S_1(x_1,...,x_n;y_1,...,y_n)$ and 
$\partial S_2(x_1,...,x_n;y_1,...,y_n)$ are $C_k$-equivalent. 

{\rm (2)} Let $v$ be a Vassiliev invariant of 
type $(k_1,...,k_l)$ for ordered, oriented links. 
Then the invariant $s_{(x_1,...,x_n;y_1,...,y_n)}$ for spatial embeddings 
of $G$ induced from $v$ is 
a Vassiliev invariant of type $(k_1,...,k_l)$.  $\Box$}

\medskip
{\bf Proof of Theorem 1.8.} Suppose that $G_k$ is abelian. 
Let $f_1$ and $f_2$ be spatial 
$\theta$-curves as illustrated in Fig. 3.5. 
Since $G_k$ is abelian, 
$g=f_1\# f_1\# f_2$ and $h=f_1\# f_2\# f_1$ 
are $C_k$-equivalent. Then, by Theorem 1.1(1), 
$g-h\in {\cal V}(\underbrace{1,...,1}_{k})$. 
This means that $g$ and $h$ cannot be distinguished by any Vassiliev 
invariant of order $\leq k-1$. Let $S_g(1,1,-1;0,0,1)$ and 
$S_h(1,1,-1;0,0,1)$ be the surfaces obtained from the 
cannonical disk/band surfaces for $g$ and $h$ respectively. 
We note that $\partial S_g(1,1,-1;0,0,1)$ and 
$\partial S_h(1,1,-1;0,0,1)$ contain 
pretzel knots $K_g=P(3,3,-3,-2)$ and 
$K_h=P(3,-3,3,-2)$ respectively, see Fig. 3.6. 
Let $v$ be a Vassiliev invariant for oriented link 
of order $\leq k-1$. 
By combining Theorem 3.9(1) and the fact that a $C_k$-moves preserves 
Vassiliev invariants of order $\leq k-1$ \cite{Habiro1} 
(or simply by Theorem 3.9(2)), we have 
$v(\partial S_g(1,1,-1;0,0,1))=v(\partial S_h(1,1,-1;0,0,1))$. 
Hence $v(P(3,3,-3,-2))=v(P(3,-3,3,-2))$. 
Let $Q_K(q)$ be the quantum invariant of a knot $K$ corresponding 
to the representation of the partition $(2,1)$ of the quantum 
enveloping algebra ${\cal U}(sl_4)$. 
Then, using the computer software \lq K2K' by 
M. Ochiai and N. Imafuji \cite{O-I}, we have 
\[\begin{array}{rcl}
\displaystyle \frac{Q_{Kg}(q)-Q_{K_h}(q)}{Q_O(q)}&=&
q^8( -1 + q)^{11}(1+q)^{11}(1+q^2)^3(1-q+q^2)^3(1+q+q^2)^3\\[-.5em]
&&\times( 1 + q^4)^2(1 - q^2 + q^4 )( 1 -q+ q^2 -q^3+ q^4 -q^5+ q^6)^2\\
&&\times( 1 +q+ q^2 +q^3+ q^4 +q^5+ q^6)^2(1+q^8)\\
&&\times(1-q^2+q^4-q^6+q^8-q^{10}+q^{12}),
\end{array}\]
where $O$ is a trivial knot. 
Since this is divisible by $(1-q)^{11}$ and 
is not divisible by $(1-q)^{12}$, these 
pretzel knots can be distinguished by a Vassiliev 
invariant of order $\leq 11$ \cite{B-L}. 
Hence we have $k< 12$. This completes the proof. $\Box$

\begin{center}
\includegraphics[trim=0mm 0mm 0mm 0mm, width=.7\linewidth]
{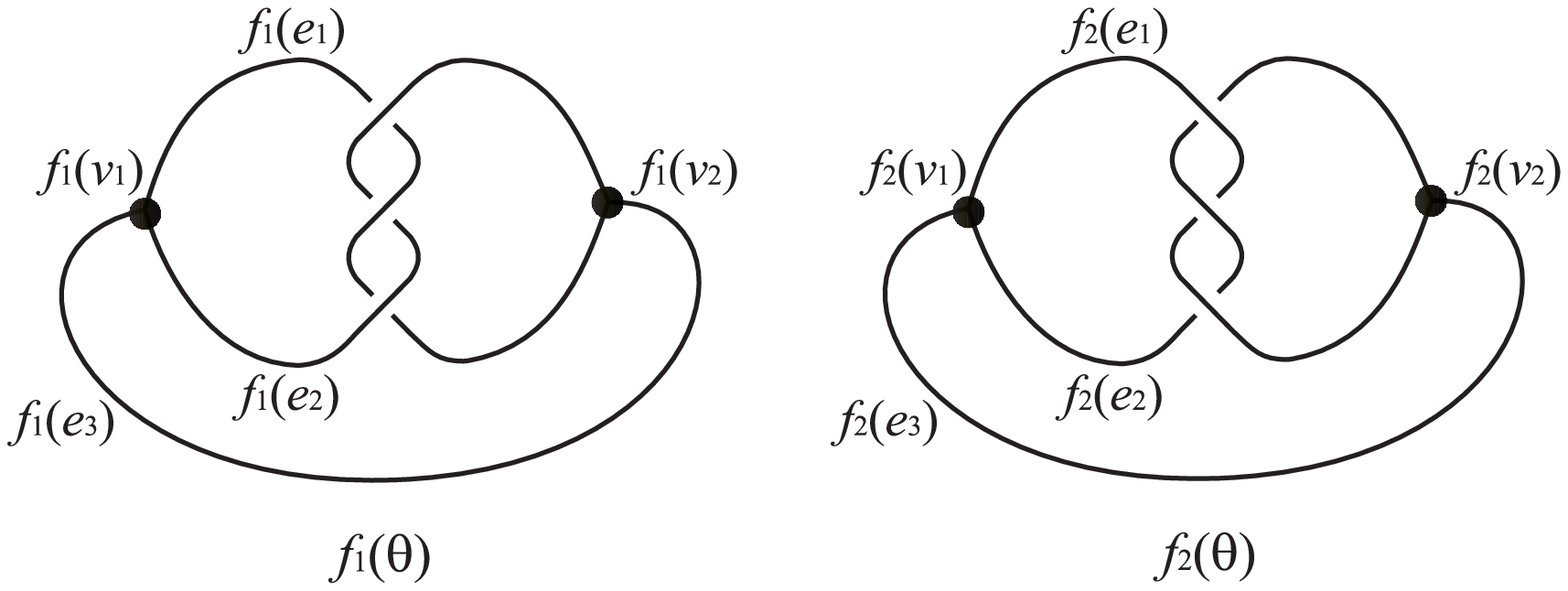}

Fig. 3.5
\end{center}
\begin{center}
\includegraphics[trim=0mm 0mm 0mm 0mm, width=.75\linewidth]
{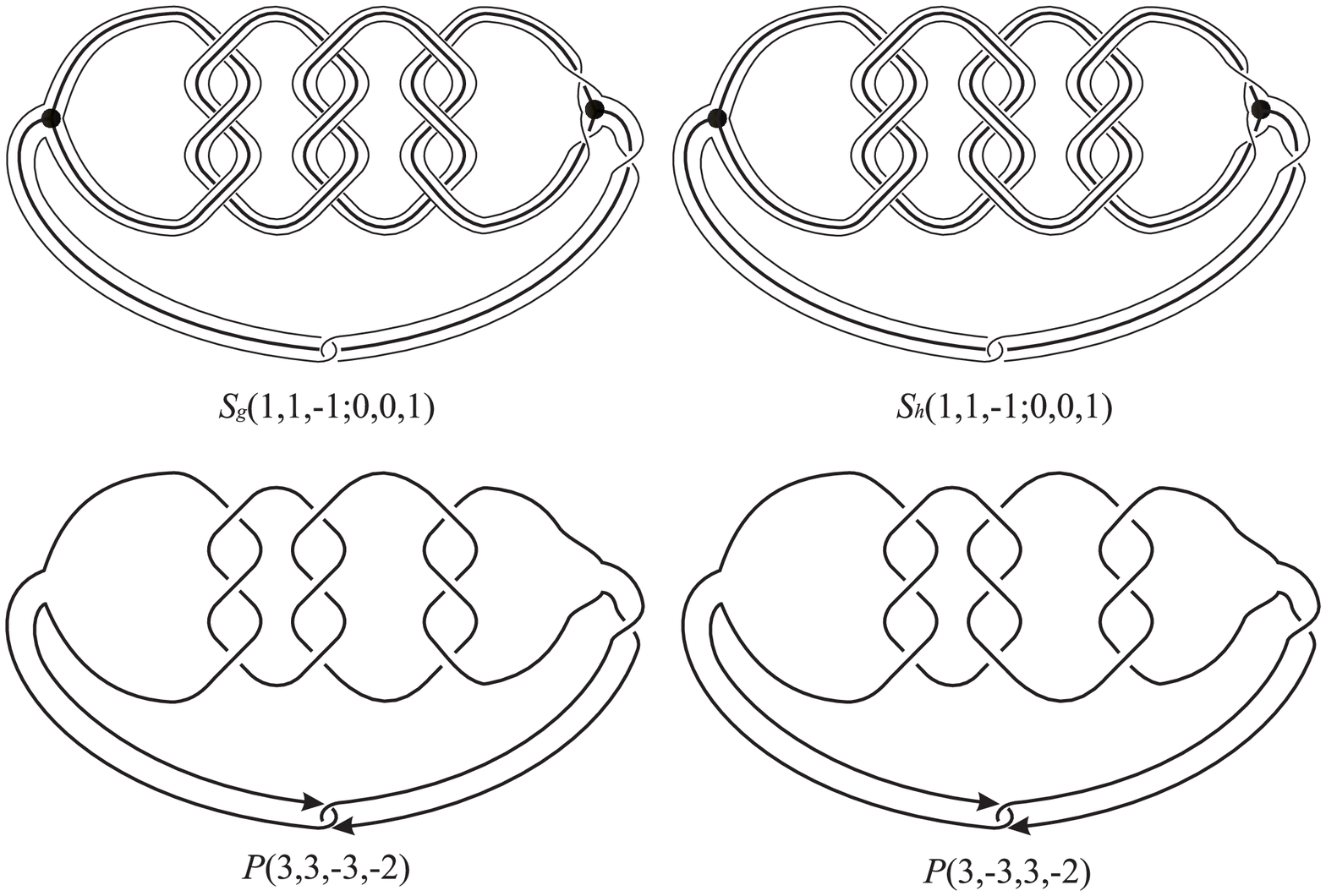}

Fig. 3.6
\end{center}

\medskip\noindent
{\bf Acknowledgement} \\
The author would like to thank 
Professor Jun Murakami, Professor J\'ozef Przytycki and 
Professor Pawe{\l} Traczyk for their valuable advice 
on distingushing two mutant knots by Vassiliev invariants. 
Their advice was useful for proving Theorem 1.8 

\bigskip
\footnotesize{
 }

\begin{thebibliography}{50}


\bibitem{B-L}J. Birman and X.-S. Lin: Knot polynomials and Vassiliev's
invariants, {\it Invent.
Math.}, 111, 225-270, 1993.



\bibitem{Habiro2}K. Habiro: Aru musubime no kyokusyo sousa no zoku ni tuite 
(in Japanese), Master thesis in Tokyo University, 1994.

\bibitem{Habiro1}K. Habiro: Claspers and finite type invariants of links,
{\it Geom. Topol.}, 4. 1-83, 2000.\\
{\bf http://www.maths.warwick.ac.uk/gt/GTVol4/paper1.abs.html}

\bibitem{Kan} T. Kanenobu: Vassiliev-type invariants of a theta-curve, 
{\it J. Knot Theory and Its Ramifications}, 6, 455-477, 1997.

\bibitem{KSWZ} L. Kauffman, J. Simon, K. Wolcott and P. Zhao: 
Invariants of theta-curves and other graphs in 3-space, 
{\it Topology Appl.} 49, 193-216, 1993.

\bibitem{Koi} A. Koike: Finite-type invariants of embeddings of a
theta-curve up to $4$, {\it Yokohama Math. J.}, 47, 245-252, 1999.


\bibitem{Mellor} B. Mellor: Finite-type link homotopy invariants II: 
Milnor's $\bar{\mu}$ invariants, preprint.\\
{\bf math.GT/9812119}



\bibitem{M-T}T. Motohashi and K. Taniyama: Delta unknotting operation 
and vertex homotopy of graphs in ${\Bbb R}^3$, 
Proceedings of Knots 96, (S. Suzuki ed.), 
World Sci. Publ. Co., 185-200, 1997.

\bibitem{O-I} M. Ochiai and N. Imafuji:  
KNOT2000(K2K)\\ 
{\bf ftp://ftp.ics.nara-wu.ac.jp/pub/ochiai}



\bibitem{SSY} T. Soma, H. Sugai and A. Yasuhara: 
Disk/band surfaces of spatial graphs, {\it Tokyo J. Math.}, 
20, 1-11, 1997.

\bibitem{Stanford0}T. Stanford: The functoriality of Vassiliev-type 
invariants of links, braids, and knotted
graphs. Random knotting and linking (Vancouver, BC, 1993). 
{\it J. Knot Theory Ramifications} 3, 247--262, 1994.

\bibitem{Stanford}T. Stanford: Finite-type invariants of knots, links and
graphs, {\it Topology}, 35, 1027-1050, 1996.



\bibitem{Stanford4}T. Stanford: Braid commutators and delta finite-type 
invariants, preprint.\\
{\bf math.GT/9907071}

\bibitem{suzuki}S. Suzuki: A prime decomposition theorem for a graph 
in $3$-sphere, Topology and Computer Science, 
(S. Suzuki ed.), Kinokuniya, 259-276, 1987. 

\bibitem{tani}K. Taniyama: Cobordism, homotopy and homology of graphs in 
${\Bbb R}^3$, {\it Topology}, 33, 509-523, 1994. 

\bibitem{T-Y0}K. Taniyama and A. Yasuhara: Local moves on spatial graphs and
finite type invariants, preprint. 

\bibitem{T-Y}K. Taniyama and A. Yasuhara: Band description of knots 
and Vassiliev invariants, to appear in 
{\it Math. Proc. Cambridge Philos. Soc.}\\
{\bf math.GT/003021}





\bibitem{Wolcott}K. Wolcott: The knotting of theta curves and other 
graphs in $S^3$, Geometry and Topology (C. McCrory and T. Shifrin ed.), 
Marcel Dekker, 325-346, 1987. 

\end{thebibliography}
\end{document}